\documentclass[11pt,leqno]{amsart}
\theoremstyle{definition}
\usepackage{indentfirst, amssymb,amsmath,amsthm}
\usepackage{csquotes}
\usepackage{hyperref}
\usepackage{fancyhdr} 
\newtheorem*{theoA}{Theorem A}
\newtheorem*{theoB}{Theorem B}

\newtheorem{theo}{Theorem}[section]
\newtheorem{lem}{Lemma}[section]
\newtheorem{cor}{Corollary}[section]

\newtheorem{defi}{Definition}[section]

\newtheorem{question}{Question}[section]
\newtheorem*{questionA}{Question A}
\newcommand{\ol}{\overline}
\newcommand{\be}{\begin{equation}}
\newcommand{\ee}{\end{equation}}
\newcommand{\beas}{\begin{eqnarray*}}
\newcommand{\eeas}{\end{eqnarray*}}
\newcommand{\bea}{\begin{eqnarray}}
\newcommand{\eea}{\end{eqnarray}}
\newcommand{\lra}{\longrightarrow}
\numberwithin{equation}{section}
\begin{document}
\title[On the cardinality of unique range sets with weight one]{On the cardinality of unique range sets with weight one}
\date{}
\author[B. Chakraborty and S. Chakraborty ]{Bikash Chakraborty$^{1}$ and Sagar Chakraborty$^{2}$ }
\date{}
\address{$^{1}$ Department of Mathematics, University of Kalyani, West Bengal, India-741 235.}
\address{$^{1}$Current Address: Department of Mathematics, Ramakrishna Mission Vivekananda Centenary College, Rahara, India-700 118.}
\address{$^{2}$Department of Education, University of Kalyani, West Bengal, India-741 235.}
\address{$^{2}$Current Address: Department of Mathematics, Jadavpur University, West Bengal, India- 700 032.}
\email{bikashchakraborty.math@yahoo.com, bikashchakrabortyy@gmail.com}
\email{sagarchakraborty55@gmail.com}
\maketitle
\let\thefootnote\relax
\begin{abstract} Two meromorphic functions $f$ and $g$ are said to share the set $S\subset \mathbb{C}\cup\{\infty\}$ with weight $l\in\mathbb{N}\cup\{0\}\cup\{\infty\}$, if $E_{f}(S,l)=E_{g}(S,l)$ where $$E_{f}(S,l)=\bigcup\limits_{a \in S}\{(z,t) \in \mathbb{C}\times\mathbb{N}~ |~ f(z)=a ~\text{with~ multiplicity}~ p\},$$
where $t=p$ if $p\leq l$ and $t=p+1$ if $p>l$.\par
In this paper, we improve and supplement the result of L. W. Liao and C. C. Yang (On the cardinality of the unique range sets for meromorphic and entire functions, Indian J. Pure appl. Math., 31 (2000), no. 4, 431-440) by showing that there exist a finite set $S$ with cardinality $\geq 13$ such that $E_{f}(S,1)=E_{g}(S,1)$ implies $f\equiv g$.
\end{abstract}
\section{Introduction and Definitions}
By $\mathbb{C}$ and $\mathbb{N}$, we mean the set of complex numbers and set of natural numbers respectively. By meromorphic function, we mean an analytic function defined on $\mathbb{C}$ except possibly at isolated singularities, each of which is a pole.\par
The tool we used in this paper is Nevanlinna Theory. For the standard notations of the Nevanlinna theory, one can go through the Hayman's Monograph (\cite{8}).\par
It will be convenient to let $E$ denote any set of positive real numbers of finite linear measure, not necessarily the same at each occurrence. For any non-constant meromorphic function $h(z)$ we denote by  $S(r,h)$ any quantity satisfying
$$S(r,h)=o(T(r,h))\;\;\;\;\;( r\lra \infty, r\not\in E).$$
Suppose $f$ and $g$ be two non-constant meromorphic functions and $a\in \mathbb{C}$. We say that $f$ and $g$ share the value $a$-CM (counting multiplicities), provided that $f-a$ and $g-a$ have the same zeros with the same multiplicities. Similarly, we say that $f$ and $g$ share the value $a$-IM (ignoring multiplicities), provided that $f-a$ and $g-a$ have the same set of zeros, where the multiplicities are not taken into account.
Moreover, we say that $f$ and $g$ share $\infty$-CM (resp. IM), if $1/f$ and $1/g$ share $0$-CM (resp. IM).\par
In the course of studying the factorization of meromorphic functions, F. Gross (\cite{7}) first generalized the idea of value sharing by introducing the concept of a unique range set. Before going to the details of the paper, we first recall the definition of set sharing.
\begin{defi}(\cite{f})
For a non-constant meromorphic function $f$ and any set $S\subset \mathbb{C}\cup\{\infty\}$, we define
$$E_{f}(S)=\bigcup\limits_{a \in S}\{(z,p) \in \mathbb{C}\times\mathbb{N}~ |~ f(z)=a ~with~ multiplicity~ p\},$$
$$\ol{E}_{f}(S)=\bigcup\limits_{a \in S}\{(z,1) \in \mathbb{C}\times\mathbb{N}~ |~ f(z)=a ~with~ multiplicity~ p\}.$$
If $E_{f}(S)=E_{g}(S)$ \big(resp. $\ol{E}_{f}(S)=\ol{E}_{g}(S)$\big), then it is said that $f$ and $g$ share the set $S$ counting multiplicities or in short CM (resp. ignoring multiplicities or in short IM).\par
Thus, if $S$ is singleton, then it coincides with the usual definition of value sharing.
\end{defi}
In 1977, F. Gross (\cite{7}) proposed the following problem which has later became popular as \enquote{\emph{Gross Problem}}. The problem was as follows:
\begin{questionA}
Does there exist a finite set $S$ such that any two non-constant entire functions $f$ and $g$ sharing the set $S$ must be $f=g$?
\end{questionA}
In 1982, F. Gross and C. C. Yang (\cite{gy}) gave the affirmative answer to the above question as follows:
\begin{theoA} (\cite{gy})  Let $S=\{z\in \mathbb{C} : e^{z}+z=0\}$. If two entire functions $f$, $g$ satisfy  $E_{f}(S)=E_{g}(S)$, then $f\equiv g$.
\end{theoA}
In that paper (\cite{gy}), they first introduced the terminology unique range set for entire function (in short URSE). Later the analogous definition for meromorphic functions was also introduced in literature.
\begin{defi}(\cite{f})
Let  $S\subset \mathbb{C}\cup\{\infty\}$ and $f$ and $g$ be two non-constant meromorphic (resp. entire) functions. If $E_{f}(S)=E_{g}(S)$ implies $f\equiv g$, then $S$ is called a unique range set for meromorphic (resp. entire) functions or in brief URSM (resp. URSE).
\end{defi}
In 1997, H. X. Yi (\cite{yi1}) introduced the analogous definition for reduced unique range sets.
\begin{defi}(\cite{yi1})
A set $S\subset \mathbb{C}\cup\{\infty\}$ is said to be a unique range set for meromorphic (resp. entire) functions in ignoring multiplicity, in short URSM-IM (resp. URSE-IM) or a reduced unique range set for meromorphic (resp. entire) functions, in short RURSM (resp. RURSE) if $\ol E_{f}(S)=\ol E_{g}(S)$ implies $f\equiv g$ for any pair of non-constant meromorphic (resp. entire) functions.
\end{defi}
During the last few years the notion of unique as well as reduced unique range sets have been generating an increasing interest among the researchers. For the literature, one can also go through the research monograph written by C. C. Yang and H. X. Yi (\cite{f}).\par Next we recall the following notations:
$$\lambda_{M}=\inf\{\sharp(S)~|~S~~\text{is~~an~~URSM}\}~~\text{and}~~\lambda_{E}=\inf\{\sharp(S)~|~S~~\text{is~~an~~URSE}\},$$
where $\sharp(S)$ is the cardinality of the set $S$.\par
In 1996, P. Li (\cite{pl}) showed that $\lambda_{E}\geq5$ whereas C. C. Yang and H. X Yi (\cite{f}) established that $\lambda_{M}\geq 6$. Also, G. Frank and M. Reinders (\cite{1}) estimated that $\lambda_{M}\leq 11$. And for entire functions corresponding estimation is $\lambda_{E}\leq 7$. Till date these estimations are the best.\par
In course of time, researchers are also paying their attention to find the lowest cardinality of URSM-IM as well as URSE-IM. In 1997, H. X. Yi (\cite{yi1}) gave the existence of URSM-IM with $19$ elements. After one year, in 1998, H. X. Yi (\cite{yi2}) further improved his result (\cite{yi1}) and obtained URSM-IM of $17$  elements. Also in 1998, M. L. Fang and H. Guo (\cite{fg}) and in 1999, S. Bartels (\cite{01}) independently gave the existence of URSM-IM with $17$ elements. In connection to our discussions, the following question is natural:
\begin{question}\label{q11} Can one further reduced the lower bound of the unique range sets by relaxing the sharing notations?
\end{question}
As an attempt to reduce the cardinalities of unique range sets, L. W. Liao and C. C. Yang (\cite{liao}) introduced the following  notation:
\begin{defi}
Let $f$ be a non-constant meromorphic function and  $S\subset \mathbb{C}\cup\{\infty\}$. We define
$$E_{1)}(S,f)=\bigcup\limits_{a \in S}E_{1)}(a,f),$$
where $E_{1)}(a,f)$ is the set of all simple $a$-points of $f$.
\end{defi}
For a positive integers $n\geq 3$ and $c\not=0,1$, we shall denote by $P(z)$, the Frank-Reinders polynomial (\cite{1}) as:
\bea\label{abcp1}
P(z)=\frac{(n-1)(n-2)}{2}z^{n}-n(n-2)z^{n-1}+\frac{n(n-1)}{2}z^{n-2}-c,
\eea
Clearly the restrictions on $c$ implies that $P(z)$ has only simple zeros. Using the methods of Frank-Reinders (\cite{1}), in connection to the Question \ref{q11},  L. W. Liao and C. C. Yang (\cite{liao}) Proved following Theorem:
\begin{theoB}(\cite{liao})
Suppose that $n(\geq1)$ be a positive integer. Further suppose that $S=\{z :P(z)=0\}$ where the polynomial $P(z)$ of degree $n$ defined by (\ref{abcp1}). Let $f$ and $g$ be two non-constant meromorphic functions satisfying $E_{1)}(S,f)=E_{1)}(S,g)$. If $n\geq15$, then $f\equiv g$.
\end{theoB}
The motivation of this paper is to improve and supplement \emph{Theorem B} utilizing the method of Frank-Reinders (\cite{1}). Before going to state our main result, we recall some well known definitions which will be useful for the proof of the main result of this paper.
\begin{defi} (\cite{lahiri})
For a non-constant meromorphic function $f$ and any set $S\subset \mathbb{C}\cup\{\infty\}$, $l\in\mathbb{N}\cup\{0\}\cup\{\infty\}$, we define
$$E_{f}(S,l)=\bigcup\limits_{a \in S}\{(z,t) \in \mathbb{C}\times\mathbb{N}~ |~ f(z)=a ~with~ multiplicity~ p\},$$
where $t=p$ if $p\leq l$ and $t=p+1$ if $p>l$.\par
Two meromorphic functions $f$ and $g$ are said to share the set $S$ with weight $l$, if $E_{f}(S,l)=E_{g}(S,l)$.\par
Clearly $E_{f}(S)=E_{f}(S,\infty)$ and $\ol{E}_{f}(S)=E_{f}(S,0)$.
\end{defi}
\begin{defi} \label{d3} Suppose $a\in\mathbb{C}\cup\{\infty\}$ and $m\in \mathbb{N}$.
\begin{enumerate}
\item [i)] We denote by $N(r,a;f\mid=1)$, the counting function of simple $a$-points of $f$,
\item [ii)] by $N(r,a;f\mid\leq m)\; (resp.~~ N(r,a;f\mid\geq m)$, we denote the counting function of those $a$-points of $f$ whose multiplicities are not greater (resp. less) than $m$ where each $a$-point is counted according to its multiplicity.
\end{enumerate}
Similarly, one can define $\ol N(r,a;f\mid\leq m)$ and $\ol N(r,a;f\mid\geq m)$ as the reduced counting function of $N(r,a;f\mid\leq m)$ and $N(r,a;f\mid\geq m)$ respectively.\par
Analogously $N(r,a;f\mid <m)$, $N(r,a;f\mid >m)$, $\ol N(r,a;f\mid <m)$ and $\ol N(r,a;f\mid >m)$ are defined.
\end{defi}
\begin{defi}\label{d5} Suppose that $f$ and $g$ be two non-constant meromorphic functions such that $f$ and $g$ share $(a,0)$. Further suppose that $z_{0}$ be an $a$-point of $f$ with multiplicity $p$, an $a$-point of $g$ with multiplicity $q$.
\begin{enumerate}
\item [i)] We denote by $\ol N_{L}(r,a;f)$, the reduced counting function of those $a$-points of $f$ and $g$ where $p>q$,
\item [ii)] by $N^{1)}_{E}(r,a;f)$, the counting function of those $a$-points of $f$ and $g$ where $p=q=1$,
\item [iii)] by $\ol N^{(2}_{E}(r,a;f)$, the reduced counting function of those $a$-points of $f$ and $g$ where $p=q\geq 2$.
\end{enumerate}
Similarly, we can define $\ol N_{L}(r,a;g)$, $N^{1)}_{E}(r,a;g)$, $\ol N^{(2}_{E}(r,a;g)$.
\end{defi}
When $f$ and $g$ share $(a,m)$, $m\geq 1$, then $N^{1)}_{E}(r,a;f)=N(r,a;f\mid=1)$.
\begin{defi} We denote by $\ol N(r,a;f\mid=k)$, the reduced counting function of those $a$-points of $f$ whose multiplicities is exactly $k$, where $k\geq 2$ is an integer.
\end{defi}
\begin{defi} \label{d7}\cite{lahiri} Let $f$, $g$ share a value $a$ IM. We denote by $\ol N_{*}(r,a;f,g)$, the reduced counting function of those $a$-points of $f$ whose multiplicities differ from the multiplicities of the corresponding $a$-points of $g$. Clearly
$$\ol N_{*}(r,a;f,g) \equiv \ol N_{*}(r,a;g,f)~~\text{and}~~\ol N_{*}(r,a;f,g)=\ol N_{L}(r,a;f)+\ol N_{L}(r,a;g).$$
\end{defi}
\section{Main Results}
\begin{theo}\label{thB3}
Suppose that $n(\geq1)$ be a positive integer. Further suppose that $S=\{z :P(z)=0\}$ where the polynomial $P(z)$ of degree $n$ defined by (\ref{abcp1}). Let $f$ and $g$ be two non-constant meromorphic functions satisfying $E_{f}(S,1)=E_{g}(S,1)$. If $n\geq13$, then $f\equiv g$.
\end{theo}
\begin{cor}\label{thB4}
Suppose that $n(\geq1)$ be a positive integer. Further suppose that $S=\{z :P(z)=0\}$ where the polynomial $P(z)$ of degree $n$ defined by (\ref{abcp1}). Let $f$ and $g$ be two non-constant entire functions satisfying $E_{f}(S,1)=E_{g}(S,1)$. If $n\geq8$, then $f\equiv g$.
\end{cor}
\section{Lemmas}
We define for any two non-constant meromorphic functions $f$ and $g$\\
 $$Q(z)=\frac{P(z)+c}{c},~ F=Q(f) ,~ G=Q(g).$$
Henceforth, we shall denote by $H$ the following three functions
$$H=\left(\frac{F^{''}}{F^{'}}-\frac{2F^{'}}{F-1}\right)-\left(\frac{G^{''}}{G^{'}}-\frac{2G^{'}}{G-1}\right).$$
\begin{lem}\label{l1}(\cite{11}) Let $f$ be a non-constant meromorphic function and let \[R(f)=\frac{\sum\limits _{k=0}^{n} a_{k}f^{k}}{\sum \limits_{j=0}^{m} b_{j}f^{j}}\] be an irreducible rational function in $f$ with constant coefficients $\{a_{k}\}$ and $\{b_{j}\}$where $a_{n}\not=0$ and $b_{m}\not=0$ Then $$T(r,R(f))=dT(r,f)+S(r,f),$$ where $d=\max\{n,m\}$.\end{lem}
\begin{lem}(\cite{f})
For a non-constant meromorphic function $f$,
$$T\left(r,\frac{1}{f}\right)=T(r,f)+O(1),$$
where $O(1)$ is a bounded quantity depending on $a$.
\end{lem}
\begin{lem}(\cite{f})
For a non-constant meromorphic function $f$ and for a complex number $a\in\mathbb{C}\cup\{\infty\}$
$$T\left(r,\frac{1}{f-a}\right)=T(r,f)+O(1),$$
where $O(1)$ is a bounded quantity depending on $a$.
\end{lem}
\begin{lem}(\cite{f}, page no. 23)\label{1}
Suppose that $f$ is a non-constant meromorphic function in the complex plane and $a_{1},a_{2},\ldots,a_{q}$ are $q(\geq3)$ distinct values in $\mathbb{C}\cup\{\infty\}$ . Then
$$(q-2)T(r,f)<\sum\limits_{j=1}^{q}\ol{N}(r,a_{j};f)+S(r,f)$$
where $S(r,f)$ is a quantity such that $\frac{S(r,f)}{T(r,f)} \to 0$ as $r \to+\infty$ out side of a set $E$ in $(0,\infty)$ with finite linear measure.
\end{lem}
\medbreak
  A polynomial $P(z)$ over $\mathbb{C}$, is called a \emph{uniqueness polynomial} for meromorphic (resp. entire) functions, if for any two non-constant meromorphic (resp. entire) functions $f$ and $g$, $P(f)\equiv P(g)$ implies $f\equiv g$.\par
In 2000, H. Fujimoto (\cite{2}) first discovered a special property of a polynomial, which was later termed as critical injection property. A polynomial $P(z)$ is said to satisfy \emph{critical injection property} if $P(\alpha )\not =P(\beta )$ for any two distinct zeros $\alpha $, $\beta $ of the derivative $P'(z)$.\par
Let $P(z)$ be a monic polynomial without multiple zero whose derivatives has mutually distinct $k$-zeros given by $d_{1}, d_{2}, \ldots, d_{k}$ with multiplicities $q_{1}, q_{2}, \ldots, q_{k}$ respectively. The following theorem of Fujimoto helps us to find many uniqueness polynomials.
\begin{lem}\label{fuji}(\cite{2a}) Suppose that $P(z)$ satisfy critical injection property.  Then $P(z)$ will be a uniqueness polynomial if and only if $$\sum \limits_{1\leq l<m\leq k}q_{_{l}}q_{m}>\sum \limits_{l=1}^{k} q_{_{l}}.$$
In particular, the above inequality is always satisfied whenever $k\geq 4$. When $k=3$ and $\max \{q_{1},q_{2},q_{3}\}\geq 2$ or when $k=2$, $\min \{q_{1},q_{2}\}\geq 2$ and $q_{1}+q_{2}\geq 5$, then also the above inequality holds.
\end{lem}
\section {Proof of the theorem}
\begin{proof} [\textbf{Proof of Theorem \ref{thB3} }]
By the assumption, it is clear that $F$ and $G$ shares $(1,1)$. Now we consider two cases.\\
\medbreak
\textbf{Case-1} First we assume that $H \not\equiv 0$. Then by simple calculations, we have
\bea\nonumber\label{cb1} N(r,\infty;H) &\leq& \ol{N}(r,0;F|\geq2)+\ol{N}(r,0;G|\geq2)+\ol{N}(r,\infty;F)\\
\nonumber&+& \ol{N}(r,\infty;G)+\ol{N}_{*}(r,1;F,G)+\ol{N}_{0}(r,0;F')+\ol{N}_{0}(r,0;G'),\eea
where $\ol{N}_{0}(r,0;F')$ is the reduced counting function of zeros of $F'$ which is not zeros of $F(F-1)$.
Thus
\bea\label{cb2} N(r,\infty;H) &\leq& \ol{N}(r,0;f)+\ol{N}(r,0;g)+\ol{N}(r,1;f)\\
\nonumber &+& \ol{N}(r,1;g)+\ol{N}(r,\infty;f)+\ol{N}(r,\infty;g)\\
\nonumber &+& \ol{N}_{*}(r,1;F,G)+\ol{N}_{\star}(r,0;f')+\ol{N}_{\star}(r,0;g'),\eea
where $\ol{N}_{\star}(r,0;f')$ is the reduced counting function of zeros of $f'$ which is not zeros of $f(f-1)$ and $(F-1)$ . We denote by $\ol N_{*}(r,1;F,G)$ the reduced counting function of those $1$-points of $F$ whose multiplicities differ from the multiplicities of the corresponding $1$-points of $G$.
Clearly \bea \label{cb2.5} \overline{N}(r,1;F|=1)=\overline{N}(r,1;G|=1)\leq N(r,\infty;H),\eea
where $\ol{N}(r,1;F|=1)$ the counting function of those simple $1$-points of $F$ which are also simple $1$-points of $G$.\par
Now using the Second Fundamental Theorem, equations (\ref{cb2}) and (\ref{cb2.5}), we get
\bea\label{n2}&& (n+1)(T(r,f)+T(r,g))\\
\nonumber &\leq& \overline{N}(r,\infty;f)+\overline{N}(r,\infty;g)+\overline{N}(r,0;f)+\overline{N}(r,0;g)\\
\nonumber &+& \ol{N}(r,1;f)+\ol{N}(r,1;g)+\overline{N}(r,1;F)+\overline{N}(r,1;G)\\
\nonumber&-& N_{\star}(r,0,f')-N_{\star}(r,0,g')+S(r,f)+S(r,g)\\
\nonumber &\leq& 2\{\overline{N}(r,\infty;f)+\overline{N}(r,\infty;g)+\overline{N}(r,0;f)+\overline{N}(r,0;g)\\
\nonumber &+& \ol{N}(r,1;f)+\ol{N}(r,1;g)\}+\{\overline{N}(r,1;F)+\overline{N}(r,1;G)\\
\nonumber &-&\overline{N}(r,1;F|=1)\}+\ol{N}_{*}(r,1;F,G)+S(r,f)+S(r,g).\eea
Again \bea\label{n2.m} && \overline{N}(r,1;F)+\overline{N}(r,1;G)-\overline{N}(r,1;F|=1)+\ol{N}_{*}(r,1;F,G)\\
\nonumber&\leq&\frac{1}{2}\{N(r,1;F)+N(r,1;G)\}+\frac{1}{2}\{\ol{N}(r,1;F|\geq 2)+\ol{N}(r,1;,G|\geq 2)\}\\
\nonumber&\leq&\frac{n}{2}\left(T(r,f)+T(r,g)\right)+\frac{1}{2}\{N(r,0;f'~|~f\not=0)+N(r,0;g'~|~g\not=0)\}.
\eea
Using (\ref{n2}) and (\ref{n2.m}), we get
\bea\label{n3} && \left(\frac{n}{2}-3\right)\left(T(r,f)+T(r,g)\right)\\
\nonumber &\leq& 2\{\overline{N}(r,\infty;f)+\overline{N}(r,\infty;g)\}+\frac{1}{2}\left\{N\left(r,0;\frac{f'}{f}\right)+N\left(r,0;\frac{g'}{g}\right)\right\}\\
\nonumber&&+S(r,f)+S(r,g)\\
\nonumber &\leq& 2\{\overline{N}(r,\infty;f)+\overline{N}(r,\infty;g)\}+\frac{1}{2}\left\{N\left(r,\infty;\frac{f'}{f}\right)+N\left(r,\infty;\frac{g'}{g}\right)\right\}\\
\nonumber&&+S(r,f)+S(r,g)\\
\nonumber &\leq& \frac{5}{2}\{\overline{N}(r,\infty;f)+\overline{N}(r,\infty;g)\}+\frac{1}{2}\{\ol{N}(r,0;f)+\ol{N}(r,0;g)\}+S(r,f)+S(r,g).
\eea
That is,
\bea\label{n3.51} && (n-7)(T(r,f)+T(r,g))\\
\nonumber &\leq& 5\{\overline{N}(r,\infty;f)+\overline{N}(r,\infty;g)\}+S(r,f)+S(r,g),
\eea
which is a contradiction if  $n\geq 13$ (resp. $8$) for meromorphic (resp. entire) case.\\
\medbreak
\textbf{Case-2} Next we assume that $H\equiv 0$. Now on  integration, two times, we have
\bea \label{pe1.1} F\equiv\frac{AG+B}{CG+D},\eea
where $A,B,C,D$ are constant satisfying $AD-BC\neq 0 $.\par
Thus applying Lemma (\ref{l1}) in equation (\ref{pe1.1}), we get
\bea\label{pe1.2} T(r,f)=T(r,g)+O(1).\eea
Next we consider the following two cases:\\
\medbreak
\textbf{Subcase-2.1} First assume $AC\neq0$. Then equation (\ref{pe1.1}) can be written as
\bea F-\frac{A}{C}=\frac{BC-AD}{C(CG+D)}.\eea
Thus $$\overline{N}\left(r,\frac{A}{C};F\right)=\overline{N}(r,\infty;G).$$
Now applying the Second Fundamental Theorem and equation (\ref{pe1.2}), we get
\beas && n T(r,f)+O(1)=T(r,F)\\
 &\leq& \overline{N}(r,\infty;F)+\overline{N}(r,0;F)+\overline{N}\left(r,\frac{A}{C};F\right)+S(r,F)\\
&\leq& \overline{N}(r,\infty;f)+\overline{N}(r,0;f)+2T(r,f)+\overline{N}(r,\infty;g)+S(r,f),\\
&\leq& 5 T(r,f)+S(r,f),\eeas
which is impossible as $n>5$.\\
\medbreak
\textbf{Subcase-2.2} Thus we consider $AC=0$. Since $AD-BC\neq0$, so $A=C=0$ never occur, thus the following two subcases are obvious:\\
\medbreak
\textbf{Subsubcase-2.2.1} First assume that $A=0$ and $C\neq0$. Then obviously $B\neq0$, hence equation (\ref{pe1.1}) can be written as
\be \label{pa1} F=\frac{1}{\gamma G+\delta},\ee
where $\gamma=\frac{C}{B}$ and $\delta=\frac{D}{B}$. Now we assume that there exist no $z_{0}$ such that $F(z_0)=1$. Then applying the Second Fundamental Theorem and equation (\ref{pe1.2}), we get
\beas &&T(r,F)\\
 &\leq& \overline{N}(r,\infty;F)+\overline{N}(r,0;F)+\overline{N}(r,1;F)+S(r,F)\\
&\leq& \overline{N}(r,\infty;f)+\overline{N}(r,0;f)+2T(r,f)+S(r,f)\\
&\leq& \frac{4}{n}T(r,F)+S(r,F),\eeas
which is impossible as $n>4$.\par
Thus there exist atleast one $z_{0}$ such that $F(z_0)=1$. Since $F$ and $G$ share $1$, hence $\gamma+\delta=1$ with $\gamma\neq0$. Thus equation (\ref{pa1}) becomes, $$F=\frac{1}{\gamma G+1-\gamma},$$
and hence $$\overline{N}(r,0;G+\frac{1-\gamma}{\gamma})=\overline{N}(r,\infty;F).$$
If $\gamma\neq1$, then the Second Fundamental Theorem and equation (\ref{pe1.2}), imply
\beas && T(r,G)\\
 &\leq& \overline{N}(r,\infty;G)+\overline{N}(r,0;G)+\overline{N}(r,0;G+\frac{1-\gamma}{\gamma})+S(r,G)\\
&\leq& \overline{N}(r,\infty;g)+\overline{N}(r,0;g)+2T(r,g)+\overline{N}(r,\infty;f)+S(r,g)\\
&\leq& \frac{5}{n}T(r,F)+S(r,F),\eeas
which is again impossible as $n>5$. Hence $\gamma=1$, i.e., $FG\equiv 1$.\par
Then \bea\label{r0} f^{n-2}\prod_{i=1}^{2}(f-\gamma_{i})~g^{n-2}\prod_{i=1}^{2}(g-\gamma_{i})\equiv\frac{4c^{2}}{(n-1)^2(n-2)^2},\eea
where $\gamma_{i}$ (i=1,2) are the roots of the equation $z^{2}-\frac{2n}{n-1}z+\frac{n}{n-2}=0$.\par
Let $z_{0}$ be a $\gamma_{i}$-point of $f$ of order $p$. Then $z_{0}$ must be a pole of $g$, (say of order $q$). Then $p=nq\geq n$.
So $$\ol{N}(r,\gamma_{i};f)\leq\frac{1}{n}N(r,\gamma_{i};f).$$
Again, let $z_{0}$ be a zero of $f$ of order $t$. Then $z_{0}$ must be a pole of $g$, (say of order $s$). Then $(n-2)t=ns$.
Thus $t>s$. Now $2s=(n-2)(t-s)\geq(n-2)$. Consequently $(n-2)t=ns$ gives $t\geq\frac{n}{2}.$
So $$\ol{N}(r,0;f)\leq\frac{2}{n}N(r,0;f).$$
Similar calculations are valid for $g$ also. Thus
\beas \ol{N}(r,\infty;f) &\leq& \ol{N}(r,0;g)+\sum\limits_{i=0}^{2}\ol{N}(r,\gamma_{i};g)\\
 &\leq& \frac{2}{n}N(r,0;g)+\frac{1}{n}\sum\limits_{i=0}^{2}N(r,\gamma_{i};g)\\
 &\leq& \frac{4}{n}T(r,g)+O(1).
\eeas
Next we apply the Second Fundamental Theorem for the identity (\ref{r0}) and we get
 \beas\label{r1} && 2T(r,f)\\
\nonumber  &\leq& \ol{N}(r,\infty;f)+\ol{N}(r,0;f)+\sum\limits_{i=0}^{2}\ol{N}(r,\gamma_{i};f)+S(r,f)\\
\nonumber  &\leq& \frac{4}{n}T(r,f)+\frac{2}{n}T(r,f)+\frac{2}{n}T(r,f)+S(r,f),
\eeas
which is a contradiction as $n\geq5$.\\
\medbreak
\textbf{Subsubcase-2.2.2} Thus we consider $A\neq0$ and $C=0$. Then obviously $D\neq0$ and equation (\ref{pe1.1}) can be written as
$$F=\lambda G+\mu,$$
where $\lambda=\frac{A}{D}$ and $\mu=\frac{B}{D}$.\par
Then obviously there exist atleast one  $z_{0}$ such that $F(z_0)=1$, otherwise we arrived at a contradiction like previous subcase. Thus $\lambda+\mu=1$ with $\lambda\neq0$. Hence \par
$$\overline{N}\left(r,0;G+\frac{1-\lambda}{\lambda}\right)=\overline{N}(r,0;F).$$
If $\lambda \neq1$, then using the Second Fundamental Theorem and equation (\ref{pe1.2}), we get
\beas && T(r,G)\\
 &\leq& \overline{N}(r,\infty;G)+\overline{N}(r,0;G)+\overline{N}\left(r,0;G+\frac{1-\lambda}{\lambda}\right)+S(r,G)\\
&\leq& \overline{N}(r,\infty;g)+\overline{N}(r,0;g)+2T(r,g)+\overline{N}(r,0;f)+2T(r,f)+S(r,g)\\
&\leq& \frac{7}{n}T(r,G)+S(r,G),\eeas
which is a contradiction as $n>7$. Thus $\lambda=1$, i.e., $F\equiv G$, i.e., $P(f)\equiv P(g)$.\par Since $P'(z)=\frac{n(n-1)(n-2)}{2}z^{n-3}(z-1)^{2}$ and $P(0)\not=P(1)$, So $P(z)$ satisfy critical injection property. Thus in view of Lemma \ref{fuji}, $P(z)$ is a uniqueness polynomial. Hence $f\equiv g$. This completes the proof of the theorem.
\end{proof}
\begin{center} {\bf Acknowledgement} \end{center}
The author is grateful to the referee for his/her valuable suggestions which considerably improved the presentation of the paper.

\end{document}